\documentclass[a4paper,12pt,reqno,english]{article}
\usepackage[utf8]{inputenc}

\usepackage{leftidx}
\usepackage{amsfonts}
\usepackage{amsmath}
\usepackage{amssymb}
\usepackage{amsthm}
\usepackage{enumitem}
\usepackage{cite}
\usepackage[top=20mm,bottom=20mm,left=20mm,right=20mm]{geometry}
\usepackage[pdfstartview=FitV,bookmarks=false]{hyperref}
\usepackage{cite}
\usepackage{tikz}
\usepackage{float}
\usetikzlibrary{braids}
\usetikzlibrary{decorations.markings}
\newmuskip\pFqmuskip

\newcommand*\pFq[6][8]{%
  \begingroup 
  \pFqmuskip=#1mu\relax
  \mathchardef\normalcomma=\mathcode`,
  \mathcode`\,=\string"8000
  \begingroup\lccode`\~=`\,
  \lowercase{\endgroup\let~}\pFqcomma
  {}_{#2}F_{#3}{\left[\genfrac..{0pt}{}{#4}{#5};#6\right]}%
  \endgroup
}
\newcommand{\pFqcomma}{{\normalcomma}\mskip\pFqmuskip}

\allowdisplaybreaks

\newcommand{\C}{\ensuremath{\mathbb{C}}}

\newcommand{\Z}{\ensuremath{\mathbb{Z}}}

\numberwithin{equation}{section}

\newtheorem{theorem}{Theorem}[section]
\newtheorem{definition}{Definition}[section]
\newtheorem{lemma}{Lemma}[section]
\newtheorem{corollary}{Corollary}[section]
\newtheorem{example}{Example}[section]
\newtheorem{remark}{Remark}[section]

\newtheorem{proposition}{Proposition}[section]

\date{}
\title {Braid monodromy for a trinomial algebraic equation by means of Mellin-Barnes integral representations  } 
\author{Mutlu KO\c{C}AR, Susumu TANAB\'E} 
\begin{document}
\maketitle
%
\begin{abstract}
In this article, we establish a braid monodromy representation of  functions satisfying an algebraic  equation containing three terms (trinomial equation).
We follow global analytic continuation of the roots to a trinomial algebraic equation that are expressed by Mellin-Barnes integral representations. 
We depict braids that arise from the monodromy around all branching points. The global braid monodromy is described in terms of rational twists of strands that yield a classical Artin braid representation of algebraic functions. As a corollary, we get a precise
description of the Galois group of a trinomial algebraic equation.
\end{abstract}
{
\center{\section{Introduction}\label{intr}}
}

In this article, we examine the braid monodromy $B_n(f)$ of 
solutions $Y_0(X), \cdots, Y_{n-1}(X)$ to a trinomial algebraic equation
\begin{equation}\label{f(X,Y)0}
f(X,Y)=Y^{n}-X^{g}Y^{p}+X^{r}=0, \quad qr-ng>0, \quad h.c.f.(p,r-g)=1
\end{equation}
with $q=n-p>1$ and $(n,p)$ coprime. At the end of the article, we consider also the case when $h.c.f. (n,p) = m>1$
(Theorem \ref{AlgBraidnoncoprime}).
We formulate the  result  expressed in terms of the braid monodromy group $B_n(f)$ consisting of  $n$ strands drawn by the roots of the equation \eqref{f(X,Y)0}.
The braid monodromy group $B_n(f)$ is a subgroup of  the Artin braid group $B_n$ of $n-$strands 
\cite{ArtinZopfe}, \cite{Birman}.

In \cite{KatoNoumi},  it has been shown that, for $(n,p)$ coprime, the global monodromy group of solutions to the algebraic equation
\begin{equation}
 \varphi(x,y) = y^n +x  y^p -1 =0
\label{KNalgeq}
\end{equation}
 admits a representation as the symmetric group ${\mathfrak S}_n.$ M.Kato and M. Noumi obtained this result in applying the monodromy representation of Pochhammer hypergeometric functions established in  \cite{BH}.
 \vspace{1.0cm}
\footnoterule
{
\footnotesize{ Key words: Algebraic functions,  braid monodromy, Pochhammer hypergeoemtric functions.}

\footnotesize{AMS Subject Classification: 14H05, 20F36, 33C20.}

{\footnotesize Acknolwedgements: ST  has been partially supported by Max Planck Institut f\"ur Mathematik,  Research fund of the Galatasaray University project number 1123 "Analysis and topology of algebraic functions and period integrals."   ST and MK have been partially supported by the Scientific and Technological Research Council of Turkey 1001 Grant No. 116F130 "Period integrals associated to algebraic varieties."}}
 
 \newpage
 Here in this article, we restrict ourselves to the case of a trinomial equation, as an equation containing $M$ terms would require the use of A-hypergeometric function in $M-2$ variables as it has been discussed in \cite{BeuAlg}. For $M=3$ the question is resolved by means of the braid monodromy of an ordinary differential equation, namely the Pochhammer hypergeometric equation.

The main contents of the article are as follows. After preliminary preparations related to Mellin-Barnes integral representations \eqref{psolns}, \eqref{qsolns} of algebraic functions,
in section \ref{main} we establish the main Theorem \ref{AlgBraid} of this note that allows us to depict the braid monodromy of all non-trivial loops.

In \cite{EsterovLangBraid}, A.Esterov and L.Lang have proven that the braid group of solutions to a degree $n$ algebraic equation with general coefficients of $\C^M,$ $M \geq 3,$ is equal to the full braid group $B_n$ under the condition of reduced support of the algebaic equation. 
In our setting \eqref{f(X,Y)0}, the condition $(n,p)$ coprime entails that its support is reduced. 
Our Corollary \ref{BnfBn} recovers Theorem 1 of 
\cite{EsterovLangBraid}. Even though we have far less loops on the moduli space of algebraic equations at our disposal than in \cite{EsterovLangBraid} where loops are chosen in $\C^3$, we conclude that the induced braid monodromy group of roots  $B_n(f)$ does recover the full braid group if $h.c.f.(n, N)=1$.

In section \ref{noncoprime}, we treat the case $h.c.f. (n,p) =m \geq 1.$ In Corollary \ref{wreathbraid}, we reestablish the result of 
\cite[Corollary 2]{EsterovLangBraid} on the wreath product of the Artin braid group. Our condition  $h.c.f. (n, r) =1$
concerns special choice of a rational curve in the coefficient moduli space $\C^3$ treated in  \cite{EsterovLangBraid}.
This means that, so far as it concerns the trinomial equation \eqref{f(X,Y)0}, our result recovers \cite[Theorem 1.5]{EsterovLangWreath} on the wreath product structure of the Galois group also. 
This consequence can be derived from the fact that our braid monodromy $B_n(f) \subset B_n$ shall fit into the exact sequence
\begin{equation}\label{exactseq}
1 \rightarrow C_n \rightarrow  B_n(f)  \overset{\phi}{\rightarrow} Gal(f) \rightarrow 1
\end{equation}
 for $C_n = Ker\; \phi$ generated by commutators and squares of the braid monodromy group $B_n(f)$, $Gal(f):$ Galois group of the equation  \eqref{f(X,Y)0}. 

 In \cite{Degt}, A.Degtyarev proposed a method to study the braid monodromy of  trigonal curves that contain the class our trinomial cubic curves \eqref{f(X,Y)0} with $n=3$. First he establishes a description of the monodromy up to its center generated by Garside braids (Example \ref{Garside}). Then he uses notions like slope and type specification of a subgroup of the Burau group $Bu_3$ \cite[Definition 2.51, 2.53]{Degt} in order to precise the monodromy.

 We also express special gratitude to Dr. Elif Segah \"Ozta\c{s} who prepared scrupulously the figures for \S \ref{main}.

\section{Algebraic functions by means of Mellin-Barnes integral representations \label{thesis}}
In this section, we shall  express the solutions of the following general trinomial algebraic equation by means of their Mellin-Barnes integral representations.

\begin{equation}\label{f(X,Y)}
f(X,Y)=Y^{n}-X^{g}Y^{p}+X^{r}=0, \quad qr-ng>0, \quad h.c.f.(p,r-g)=1,
\end{equation}
with $q=n-p>1$ and $(n,p)$ coprime.

In order to obtain such an expression, we will make use of change of variables and Mellin inversion formula established in \cite{Mellin15}, \cite{Mellin21}, \cite{Bel}, \cite{STsMHHG}.

It is well known that the local behaviour of roots of \eqref{f(X,Y)} can be described by the Newton polygon method 
 \cite[8.3]{Bris.}.

Let us set the positive integer
\begin{equation}\label{N} 
qr-ng=N.
\end{equation}
 We remark here that the condition $N>0$ is required to ensure the convexity of the Newton boundary of the polynomial $f(X,Y)$. in fact, the area of the Newton polygon $$\Delta(f) = c.h. \{(0,n), (g,p), (r,0)\}$$
equals to $\frac{N}{2!}.$  See Figure \ref{figN(f(X,Y))}.

To consider the equation \eqref{f(X,Y)} with $qr-ng <0$, it is enough to switch the role
of $(X,Y) = (0,0)$ and $(X,Y) = (\infty,\infty)$ with the aid of a change of variables $(X,Y) \rightarrow (\frac{1}{X}, \frac{1}{Y}).$

In general, we impose also the codition
\begin{equation}\label{gcdnN}
h.c.f. (n, N)=h.c.f. (n,r)= 1
\end{equation}
under which the series \eqref{Capital} give a complete set of roots of the equation \eqref{f(X,Y)}.

By a suitable change of variables, we can convert the equation $\varphi(x,y)=0$ \eqref{KNalgeq} into the equation $f(X,Y)=0,$ \eqref{f(X,Y)}.

For the following homotheties 

\begin{equation}\notag
y=\lambda Y \quad and \quad x=\mu X, \quad (\lambda,\mu \neq 0)
\end{equation}
we have \begin{equation*}
    Y^{n}+\lambda^{p-n}\mu^g X^{g}Y^{p}-\lambda^{-n}=Y^{n}-X^{g}Y^{p}+X^{r}.
\end{equation*}
Thus we get the following equalities
\begin{equation}\label{eqn 4.3.10}
Y(X)=(-1)^{\frac{1}{n}}X^{\frac{r}{n}}y(x)
\end{equation}
and \begin{equation}\label{eqn 4.3.11}
x=(-1)^{\frac{p}{n}}X^{-\frac{N}{n}}.
\end{equation}

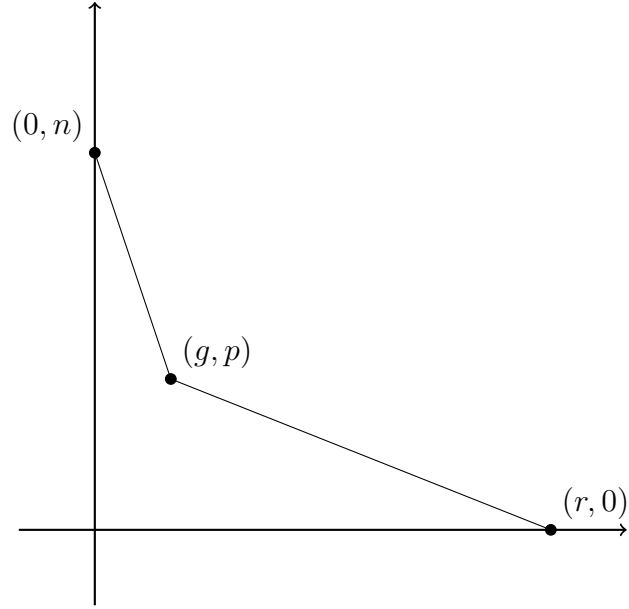
\begin{figure}
\begin{center}
\begin{tikzpicture}
\draw [thin,color=gray,step=1cm](-2,-2)  (6,6);
\draw [thick,->] (-1,0)--(7,0);
\draw [thick,->] (0,-1)--(0,7);
\draw  (0,5) node[above left]{$(0,n)$};
\draw  (1,2) node[above right]{$(g,p)$};
\draw  (6,0) node[above right]{$(r,0)$};
\draw[fill=black] (6,0) circle (2pt);
\draw[fill=black] (0,5) circle (2pt);
\draw[fill=black] (1,2) circle (2pt);
\draw (1,7/2) node[above]{};
\draw (7/2,3/2)node[above]{};
\coordinate (A) at (0,5);
\coordinate (B) at (1,2);
\coordinate (C) at (6,0);

\draw[black,thin] (A)--(B);
\draw[black,thin] (B)--(C);
\end{tikzpicture}
\end{center}
\caption{Newton polygon of $f(X,Y)$}\label{figN(f(X,Y))}
\end{figure}

In order to obtain the general hypergeometric differential equation satisfied by $Y(X)$, first we recall the following equation 
\begin{equation}
\left[\prod_{l=0}^{n-1}(\vartheta_x-l)-\frac{(-p)^{p}q^{q}}{n^n}x^{n} \prod_{j=0}^{p-1}(\vartheta_x+\frac{1+nj}{p}) \prod_{k=0}^{q-1}(\vartheta_x+\frac{-1+nk}{q})\right]y(x)=0,
\end{equation}
established in \cite{Mellin15}, \cite{Mellin21},  \cite[ Chapitre IV, 19]{Bel}, \cite[Proposition 2.6]{KatoNoumi}. By using \eqref{eqn 4.3.11} we may write 
 \begin{equation*}
   x \frac{\partial }{\partial x} = \vartheta_x=-\frac{n}{N}\vartheta_X,
\end{equation*}
and hence
\begin{equation}
\left[\frac{X^{N}}{R}\prod_{l=0}^{n-1}(\frac{n}{N}\vartheta_X+l)
- \prod_{j=0}^{p-1}(\frac{n}{N}\vartheta_X-\frac{1+nj}{p}) \prod_{k=0}^{q-1}(\frac{n}{N}\vartheta_X-\frac{-1+nk}{q})\right]y((-1)^{\frac{p}{n}}X^{-\frac{N}{n}})=0,\label{diff eqn infinty}
\end{equation}
for \begin{equation}\label{R}
    R=\frac{p^p q^q}{n^n}.
\end{equation}

Applying the change of variable $R^{-1}X^{N}=s$, hence $\vartheta_X=N \vartheta_s$ yields the following hypergeometric differential equation 
\begin{equation}\label{HGEqn-1}
\left[s\prod_{l=0}^{n-1}(\vartheta_{s}+\frac{l}{n})- \prod_{j=0}^{p-1}(\vartheta_{s}-\frac{1}{np}-\frac{j}{p}) \prod_{k=0}^{q-1}(\vartheta_{s}-(-\frac{1}{nq}+\frac{k}{q})\right]y=0.
\end{equation} 
For a further use, we remark that \begin{equation*}
    Y(X)=(-X^r)^{\frac{1}{n}}y((-1)^{\frac{p}{n}}X^{-\frac{N}{n}})
\end{equation*}
satisfies the following equation

\begin{equation}
\left[\frac{X^{N}}{R}\prod_{l=0}^{n-1}(\frac{n\vartheta_{X}-r}{N}+l)
- \prod_{j=0}^{p-1}(\frac{n\vartheta_{X}-r}{N}-\frac{1+nj}{p}) \prod_{k=0}^{q-1}(\frac{n\vartheta_{X}-r}{N}-\frac{-1+nk}{q})\right]Y(X)=0.
\end{equation}
In other words,

\begin{equation}
\left[s\prod_{l=0}^{n-1}(\vartheta_{s}+\frac{l}{n}-\frac{1}{nq}) -\prod_{j=0}^{p-1}(\vartheta_{s}-\frac{1}{pq}-\frac{j}{p}) \prod_{k=0}^{q-1}(\vartheta_{s}-\frac{k}{q})\right]s^{-\frac{ n g}{q N}}Y((Rs)^{\frac{1}{N}})=0.
\end{equation} 
This equation is reducible and we get the following equation of order $(n-1):$
\begin{equation}\label{diffeqn}
\left[s\prod_{l=0,\neq \tilde{l}}^{n-1}(\vartheta_{s}+\frac{l}{n}-\frac{1}{nq}) -\prod_{j=0}^{p-1}(\vartheta_{s}-\frac{1}{pq}-\frac{j}{p}) \prod_{k=0, \neq \tilde{r}}^{q-1}(\vartheta_{s}-\frac{k}{q})\right]s^{-\frac{ n g}{q N}}Y((R s)^{\frac{1}{N}})=0,
\end{equation} 
for an unique pair $\tilde{r}\in [0;q-1]$, $\tilde{l}\in [0;n-1]$ satisfying $n\tilde{r}+q\tilde{l}=nq+1$ as $(n,q)$ coprime.
From here on, we use the notation $[a;b] =\{a, a-1, \cdots, b-1, b\}$ for two integers $a \leq b.$

To obtain a series expansion of $Y(X)$ in $|\frac{X^{N}}{R}|>1$, we recall its Mellin-Barnes integral representation obtained in \cite{Kocar} that is based on \cite{Mellin15}, \cite{Mellin21}, \cite[Chapitre IV]{Bel} :

\begin{equation}\label{YXintegral}
    Y(X)=\frac{(-X^r)^{\frac{1}{n}}}{2\pi i n}\int_{c-i \infty}^{c+i \infty} \frac{\Gamma(z)\Gamma(\frac{1-pz}{n})}{\Gamma(\frac{1+qz}{n}+1)}e^{-\frac{pz\pi i}{n}}X^{\frac{N z}{n}}dz.
\end{equation}
On summing up residues at $z=-k$, $k\in \mathbb{Z}_{\geq 0},$ we get a series expression of one of the roots \eqref{eqn 4.3.10} near $X=\infty$:
\begin{equation}
    Y_{0}(X)=\frac{(-1)^\frac{1}{n}}{n}\sum_{k=0}^{\infty}\frac{\Gamma(\frac{1+pk}{n})(-1)^{-\frac{qk}{n}}}{\Gamma(\frac{1-qk}{n}+1)k!}(X^{\frac{1}{n}})^{-N k+r}
\end{equation}

Also we get the series \begin{equation}\label{Capital}
    Y_t(X):=Y_{0}(e(t)X),
\end{equation}
or, for $c\in (0,\frac{1}{p}), $
\begin{equation}\label{YtSeries}
    Y_{t}(X)=\frac{(-1)^\frac{1}{n}}{2\pi i n}\int_{c-i \infty}^{c+i \infty}\frac{\Gamma(z)\Gamma(\frac{1-pz}{n})}{\Gamma(\frac{1+qz}{n}+1)}e\left((N t-\frac{p}{2})\frac{z}{n}+ \frac{rt}{n}\right)X^{\frac{N z+r}{n}}dz,
\end{equation}
that is a result of a local monodromy action $X^{\frac{1}{n}}\mapsto e(\frac{t}{n})X^{\frac{1}{n}}$ for $t\in [0;p-1]$. From here on we use the notation \begin{equation}
    e(\alpha)=e^{2\pi i \alpha}.
\end{equation}
As a result we get the Newton-Puiseux series near $X=0$ of the root ${Y}_{t}(X), t\in [0;p-1]$ of an algebraic equation $f(X,Y)=0,$ 
\eqref{f(X,Y)}.

To define the Mellin-Barnes integral representation of  ${Y}_{t'}(X), t' \in [p;n-1],$
  we multiply \eqref{YtSeries} by the following meromorphic function satisfying $\phi(z+1) = \phi(z),$
\begin{equation}\label{periodicfactor}
    \phi(z) =e(\frac{z}{2})\frac{\sin \pi (\frac{(1-pz)}{n})}{ \sin \pi (\frac{-(1+qz)}{n})}.   
\end{equation}
to get 
\begin{equation}\label{qYXintegral}
    Y_{t'}(X)=\frac{\left(- ( e(t')X )^r \right)^{\frac{1}{n}}}{2\pi i n}\int_{c-i \infty}^{c+i \infty} \frac{\Gamma(z)\Gamma(\frac{-1-qz}{n})}{\Gamma(\frac{-1+pz}{n}+1)}e(\left( \frac{q}{2} + N t'\right)\frac{z}{n})X^{\frac{N z}{n}}dz.
\end{equation}
with $t' \in [p;n-1].$ The summation of residues for \eqref{qYXintegral} at $z=-k$, $k\in \mathbb{Z}_{\geq 0},$
  also gives us the series expansion \eqref{Capital}
\begin{equation}\label{qCapital}
 Y_{t'}(X) = Y_0(e(t')X)
 \end{equation}
with $t' \in [p; n-1].$ The unified representation \eqref{Capital}, \label{qCapital} valid for all roots $Y_{t}(X), t \in [0;n-1]$ based on
$Y_0(X)$ is explained by the choice of 
 a proper periodic function $\phi(z)$ \eqref{periodicfactor}.


Due to the use of two different Mellin-Barnes integrals \eqref{YtSeries}, \eqref{qYXintegral}, we shall treat
 two types of analytic expansions near $X=0$.

To achieve the analytic continuation of $p-$solutions \eqref{YtSeries} and $q-$solutions \eqref{qYXintegral}
from $X=\infty$ to $X=0,$ we use the technique of Mellin-Barnes contour throw. We get the series expansion near $X=0$ that is
obtained as an analytic continuation
of  \eqref{Capital} by summing up the residues of  \eqref{YtSeries} and \eqref{qYXintegral} at poles with positive real parts. 
See \cite[Section 3 "Solutions at the Origin and the Point at Infinity",{\bf 16-18}]{Norlund}.

 If we gather the residues of \eqref{YtSeries} at $z= \frac{1+nk}{p}, k \in \Z_{\geq 0},$ the following series expansions convergent in the disc $|\frac{X^N}{R}| <1$  (named $p-$solutions) is obtained,
\begin{equation}\label{psolns}
    Y_{t}(X)=\frac{X^\frac{r-g}{p}}{p}\sum_{k=0}^{\infty}\frac{\Gamma(\frac{1+nk}{p}) }{\Gamma(\frac{1+p+qk}{p})k!}e(\frac{ (r-g+Nk)t }{p})(X^{\frac{1}{p}})^{N k}
\end{equation}
with $t \in [0;p-1].$

If we sum up the residues of \eqref{qYXintegral} at $z= \frac{-1+nk}{q}, k \in \Z_{\geq 0},$
 we get another series expansions (named $q-$solutions) convergent in $|\frac{X^N}{R}| <1$,
\begin{equation}\label{qsolns}
    Y_{t'}(X)=\frac{X^{\frac{g}{q}}}{q}\sum_{k=0}^{\infty}\frac{\Gamma(\frac{-1+nk}{q})}{\Gamma(\frac{-1+pk}{q}+1)k!}e(\frac{(g + Nk) t'}{q})(X^{\frac{1}{q}})^{N k}.
\end{equation}
The initial exponent $\frac{r-g}{p}$ of \eqref{psolns} (resp. $\frac{g}{q}$) corresponds to the solope
of an edge of the Newton polygon of $f(X,Y)$ (Figure \ref{figN(f(X,Y))}). This matches the Newton--Puiseux expansion
near $X=0$ established by the 
classical method \cite[8.3]{Bris.}.

\section{Braids depicted by roots of a trinomial equation \label{main}}
In this section, we establish the main result of this note, namely we pursue the braids depicted by the roots $Y_t(X), t \in [0;n-1] $ as $X$ avoids their branching points
\begin{equation}
    \Sigma \cup \{0, \infty\}=\{X\in \mathbb{P}^{1}  ; X=0,\omega_0,\cdots,\omega_{N -1}, \infty\}, 
\end{equation}
where
\begin{equation}\label{branching points}
 \Sigma =\{   \omega_\ell=R^{\frac{1}{N}}e({\frac{\ell}{N}}), \quad \ell \in [0;N-1] \}.
\end{equation}
At each branching point, there exists a pair of $(t,t') \in [0;n-1]^2$ such that 
$Y_t(\omega_\ell) = Y_{t'}(\omega_\ell).$ No other roots concidence happens at $X= \omega_\ell.$
This can be seen from the fact that the system of equations  $f(X,Y)=  \frac{\partial}{\partial Y} f(X,Y) = \frac{\partial^2}{\partial Y^2} f(X,Y) =0$ has no solution for \eqref{f(X,Y)}.

 By \eqref{YtSeries}, \eqref{qYXintegral},  we get the following series expression for the algebraic function in question.
\begin{lemma}\label{Ytinfinity}
We have in $|\frac{X^{N}}{R}|>1$
 \begin{equation}\label{Ytinfty}
    Y_{t}(X)=\frac{(-1)^{\frac{1}{n}}e(\frac{rt}{n})}{n}X^{\frac{r}{n}}\sum_{k=0}^{\infty}\frac{\Gamma(\frac{1+pk}{n})}{\Gamma(\frac{1-qk}{n}+1)}(e^{\frac{-(q+2t N)\pi i}{n}}X^{-\frac{N}{n}})^{k}.
\end{equation}
\end{lemma}

Now we ask the following question. For $t \neq t^{\prime}$, which pairs of  $Y_t(X)=Y_{t^{\prime}}(X)$ coincide? That is, at $X=\omega_\ell$ which analytic continuation of algebraic solutions $Y_t(X) $ and $Y_{t^{\prime}}(X)$ \eqref{Capital} encounter to form a nontrivial
braid?

In order to provide an explicit answer to the question above, we follow the ideas and notations in \cite{KatoNoumi}.
\newline
After \cite[Proposition 2.2]{ KatoNoumi}, we consider a generalized binomial series \begin{equation} \label{psi}
    \psi(\alpha,s,x)=\sum_{k=0}^{\infty}\frac{\alpha \Gamma(\alpha + k(s+1))}{\Gamma(\alpha+ks+1)}\frac{x^k}{k!},
\end{equation}
for any complex numbers $ \alpha,s$ with the radius of convergence $\left|\frac{s^s}{(s+1)^{s+1}}\right|$ and $s,s+1 \neq 0$.
\newline
Let us denote by $y_{j}(x)$, the algebraic functions satisfying the equation \eqref{KNalgeq}, i.e, \begin{equation*}
    (y_j(x))^n+x(y_j(x))^p-1=0,
\end{equation*}
for $j \in [0;n-1]$.
\newline
The algebraic functions $y_j(x)$ admit an expression in terms of the generalized binomial function $\psi$ \cite[Proposition 2.6]{KatoNoumi} like
\begin{equation}\label{binomialpsi}
    y_j(x)= e\left(\frac{j}{n}\right) \psi\left(-\frac{1}{n},-\frac{p}{n},e\left(\frac{pj}{n}\right)x\right).
\end{equation}
\begin{lemma}\label{psirelation}
    The generalized binomial function $\psi$ satisfies the following equality \begin{equation}
        e\left(-\frac{1}{2n}\right) \psi\left(-\frac{1}{n},-\frac{p}{n},e\left(-\frac{p}{2n}\right)R^{-\frac{1}{n}}\right)=e\left(\frac{1}{2n}\right)\psi\left(-\frac{1}{n},-\frac{p}{n},e\left(\frac{p}{2n}\right)R^{-\frac{1}{n}}\right).
    \end{equation}
    Moreover, the values of both LHS and RHS are real.
\end{lemma} 
\begin{proof}
    By  \cite[Lemma 3.3]{KatoNoumi}, we know that the algebraic functions $y_j(x)$ and $y_{j+1}(x)$ coincide at the branching points $x_j=e\left(\frac{-p(1+2j)}{2n}\right)R^{-\frac{1}{n}}$ of the equation \eqref{KNalgeq}. That is, \begin{equation}\label{algbranchingpoints}
        y_j(x_j)= e\left(\frac{j}{n}\right) \psi\left(-\frac{1}{n},-\frac{p}{n},e\left(\frac{pj}{n}\right)x_j\right)= e\left(\frac{j+1}{n}\right) \psi\left(-\frac{1}{n},-\frac{p}{n},e\left(\frac{p(j+1)}{n}\right)x_j\right)=y_{j+1}(x_j).
    \end{equation}
    Here we remark that the numbering $\{x_j \}_{j=0}^{n-1}$ of branching points gives \begin{equation}
        Arg(x_{j+1})-Arg(x_j)=-  2\pi \left(\frac{p}{n}\right).
    \end{equation}
    Now the first claim follows  from a replacement $x_j=e\left(\frac{-p(1+2j)}{2n}\right)R^{-\frac{1}{n}}$ in \eqref{algbranchingpoints}. 
    By definition of $\psi$, we can write 
    $$
      e\left(-\frac{1}{2n}\right)\psi\left(-\frac{1}{n},-\frac{p}{n},e\left(-\frac{p}{2n}\right)R^{-\frac{1}{n}}\right)
      =\sum_{k=0}^{\infty}\frac{\Gamma(\frac{-1+qk}{n})e\left(-\frac{(1+pk)}{2n}\right)}{\Gamma(-\frac{1+pk}{n}+1)}\frac{R^{-\frac{k}{n}}}{k!}$$
      $$  = e\left(\frac{1}{2n}\right) \psi\left(-\frac{1}{n},-\frac{p}{n},e\left(\frac{p}{2n}\right)R^{-\frac{1}{n}}\right)
      =\sum_{k=0}^{\infty}\frac{\Gamma(\frac{-1+qk}{n})e\left(\frac{1+pk}{2n}\right)}{\Gamma(-\frac{1+pk}{n}+1)}\frac{R^{-\frac{k}{n}}}{k!}$$
    and we see that the above series are equal as complex conjugates. Therefore, they must be both real valued.
\end{proof}

A comparison of series  \eqref{Ytinfinity} and \eqref{binomialpsi}, together with the change of variables \eqref{eqn 4.3.10}, \eqref{eqn 4.3.11}, gives rise to an expression 
\begin{equation}\label{bigypsi}
    Y_{t}(X) =e\left(\frac{1+2rt}{2n}\right)X^{\frac{r}{n}} \psi\left(-\frac{1}{n},-\frac{p}{n},e \left(\frac{p(1+2rt )}{2n}\right)X^{-\frac{N}{n}}\right), \quad rt \in[0;n-1],
\end{equation}
where $Y_{t}(X)$ on the LHS is defined after \eqref{Capital}.

Due to the condition that the pair $(n,r)$ is coprime \eqref{N},
the multipication $t\rightarrow rt =j$ gives a bijective mapping from $\Z/n\Z$ to itself.

Now we arrange the branching points \eqref{branching points} as follows: 
\begin{equation}\label{rearrangedbranching}
    X_\ell^{(j)}=R^{\frac{1}{N}}e\left(\frac{p(j+1)+\ell}{N}\right)=\omega_{p(j+1)+\ell},
\end{equation}
where $j\in [0;n-1], \ell \in [0;p-1]$.
Note that for every fixed $j \in [0;n-1] $ we have the equality of sets of branching points $\{X_\ell^{(j)}\}_{\ell=0}^{N-1}=\{X_\ell^{(0)}\}_{\ell=0}^{N-1}$ with $X_\ell^{(j)}=X_{\ell+pj}^{(0)}$. 

By \eqref{bigypsi} and \eqref{rearrangedbranching}, we have \begin{equation}
    Y_{t}(X_\ell^{(j)})=e\left(\frac{1+2 rt }{2n}+\frac{(p(j+1)+\ell)r}{n N}\right)R^{\frac{r}{n N}} \psi\left(-\frac{1}{n},-\frac{p}{n},e\left(\frac{p}{2n}+\frac{p(rt-j-1)-\ell}{n}\right)R^{-\frac{1}{n}}\right).
\end{equation}
For fixed $j$ and $\ell$, we are now in a position to look for an index pair $({ t},{ t}^{\prime}) \in [0;n-1]$ such that $Y_{t}(X_\ell^{(j)})=Y_{{ t}^{\prime}}(X_\ell^{(j)})$.
\begin{proposition}\label{Ytcoincidence}
    The coincidence $Y_{t}(X_{\ell}^{(j)})=Y_{{t}^{\prime}}(X_{\ell}^{(j)})$ of roots at the branching point occurs only 
    for $({ t},{ t}^{\prime})$ satisfying $$r t^{\prime}\equiv r{ t}+1\; \; mod \; n$$ and 
    \begin{equation}\label{ptpj}
        pr { t}  \equiv p j+\ell \mod n.
    \end{equation}
    That means
    \begin{equation}\label{YtYt+1coincidence}
        Y_{t}(\omega_{\bar\ell})=Y_{{ t'}}(\omega_{\bar\ell}) \quad \text{with} \quad \bar\ell=p(j+1)+\ell \equiv p r {t'} \quad \mod n
    \end{equation}
    where $\omega_{\bar\ell}=\omega_{p(j+1)+\ell}=X_\ell^{(j)}$  for $\ell \in [0; p-1]$.
\end{proposition}

\begin{proof}
   Writing down the series expressions of $Y_{ t}(X_\ell^{(j)})$ and $Y_{{t}^{\prime}}(X_\ell^{(j)}),$ we get the following equalities 
   \begin{align*}
       Y_{t}(X_\ell^{(j)})&=e\left(\frac{1+2rt}{2n}\right)\sum_{k=0}^{\infty}\frac{\Gamma(\frac{-1+qk}{n})}{\Gamma(-\frac{1+pk}{n}+1)}e\left(\left(\frac{p}{2n}+\frac{p(rt-j-1)-\ell}{n}\right)k\right)\frac{R^{-\frac{k}{N}}}{k!}\\
       &=e\left(\frac{1+2rt^{\prime}}{2n}\right)\sum_{k=0}^{\infty}\frac{\Gamma(\frac{-1+qk}{n})}{\Gamma(-\frac{1+pk}{n}+1)}e\left(\left(\frac{p}{2n}+\frac{p(rt^{\prime}-j-1)-\ell}{n}\right)k\right)\frac{R^{-\frac{k}{N}}}{k!}\\
       &=Y_{t^{\prime}}(X_\ell^{(j)}).
   \end{align*}
   Equivalently, \begin{align*}
       \sum_{k=0}^{\infty}\frac{\Gamma(\frac{-1+qk}{n})}{\Gamma(-\frac{1+pk}{n}+1)}&e\left(\left(\frac{p}{2n}+\frac{p(rt-j-1)-\ell}{n}\right)k\right)\frac{R^{-\frac{k}{N}}}{k!}
       \\
       &=e\left(\frac{r(t^{\prime}-t) }{n}\right) \sum_{k=0}^{\infty}\frac{\Gamma(\frac{-1+qk}{n})}{\Gamma(-\frac{1+pk}{n}+1)}e\left(\left(\frac{p}{2n}+\frac{p(rt^{\prime}-j-1)-\ell}{n}\right)k\right)\frac{R^{-\frac{k}{N}}}{k!}.
   \end{align*}
       Now comparing the latter equality with the expression in Lemma \ref{psirelation} we see the following sufficient condition 
        \begin{equation*}
           p+2p(rt-j-1)-2\ell \equiv -p \mod 2n,
       \end{equation*}
    \begin{equation*}
           p+2p(rt^{\prime}-j-1)-2\ell \equiv p \mod 2n,
       \end{equation*}
       \begin{equation*}
          rt^{\prime}-rt \equiv 1 \mod n. 
       \end{equation*}
      On combining the above congruence relations, we obtain \begin{equation*}
        Y_{t}(\omega_{\bar\ell})=Y_{t'}(\omega_{\bar\ell}) \quad \text{with} \quad \bar\ell=p(j+1)+\ell \equiv prt' \mod n.  \end{equation*}
        Thus the condition \eqref{ptpj} is sufficient for the coincidence of roots.
        
        For every branching point ${ \omega}_{\bar \ell}$ there exists a unique pair $({ t}, { t'})$
        such that   $Y_{t}(\omega_{\bar\ell})=Y_{t'}(\omega_{\bar\ell})$ hold as it has been remarked after \eqref{branching points}. The pair  $({ t}, { t'})$
        found above satisfies this condition. 
       \end{proof}

\begin{remark}\label{isotopy}
Further on, braids drawn by $\{Y_t(X)\}_{t\in [0; n-1]}$
shall be understood up to braid isotopy. To a continuous path ${\bf \lambda} = \{ \lambda(s); s \in [0,1]\}$ on the $\C_X$
plane, we associate $n-$strings $ {\bf \mathbb Y (\lambda) } = \{(s, Y_t(\lambda(s)))_{t\in [0; n-1]}; s \in [0,1] \}$ in $[0,1] \times \C_Y.$
We look at another continuous $n-$strings given in the form $ {\bf \Lambda }= \{ (s, \Lambda_t (s) )_{t\in [0; n-1]} ; s \in [0,1] \}$ in $[0,1] \times \C_Y$
such that $ \Lambda_t(0) = Y_t(\lambda(0)),$ $ \Lambda_t(1) = Y_t(\lambda(1)),$ $\forall t \in [0;n-1].$ Two braids
$ {\bf \mathbb Y (\lambda) } $  and 
$ {\bf \Lambda } $ are isotopic if $ ((0,1) \times \C_Y) \setminus {\bf \mathbb Y (\lambda) }$
is homeomorphic to $ ((0,1) \times \C_Y) \setminus {\bf \Lambda }.$

See a similar notion of the equivalence of annular braids $CB_n$ in \cite[page 86]{KentPeifer}.
\end{remark}

Proposition \ref{Ytcoincidence} entails the following lemma in view of the fact that the algebraic functions $Y_t(X), t \in [0;n-1]$
are holomorphic on a domain with no intersection with $\Sigma \cup \{0\}.$

\begin{lemma}\label{Ytinterchange}
    As $X$ turns around one of the branching points $\omega_{\bar \ell}= \omega_{p(j+1)+\ell} \in \Sigma$ in anticlockwise way, the values $Y_{t}(X)$ and $Y_{t'}(X)$ interchange their relative positions for $(t, t')$  determined by the relations \eqref{ptpj},  \eqref{YtYt+1coincidence}. In the course of this position change, both branches
make an anticlockwise move with respect to another branch. The relative positions of branches $\{Y_j(X)\}_{j \in [0;n-1] \setminus \{{ t},{t'} \}}$ do not change as $X$ turns around  $\omega_{\bar \ell}$ in anticlockwise way. The movement of roots shall be understood up to braid isotopy in the sense of Remark \ref{isotopy}.
\end{lemma}

As for the formulation of the above Lemma, see the remark made after \eqref{branching points}.

In \cite[3.2 Tropical construction of coarse braids]{EsterovLangBraid} an analogous result is obtained by a minute following of
argument change of algebraic functions together with their logarithmic image contained in an amoeba of a  straight line.

\begin{figure}[ht!]
\centering
\includegraphics[scale=0.4]{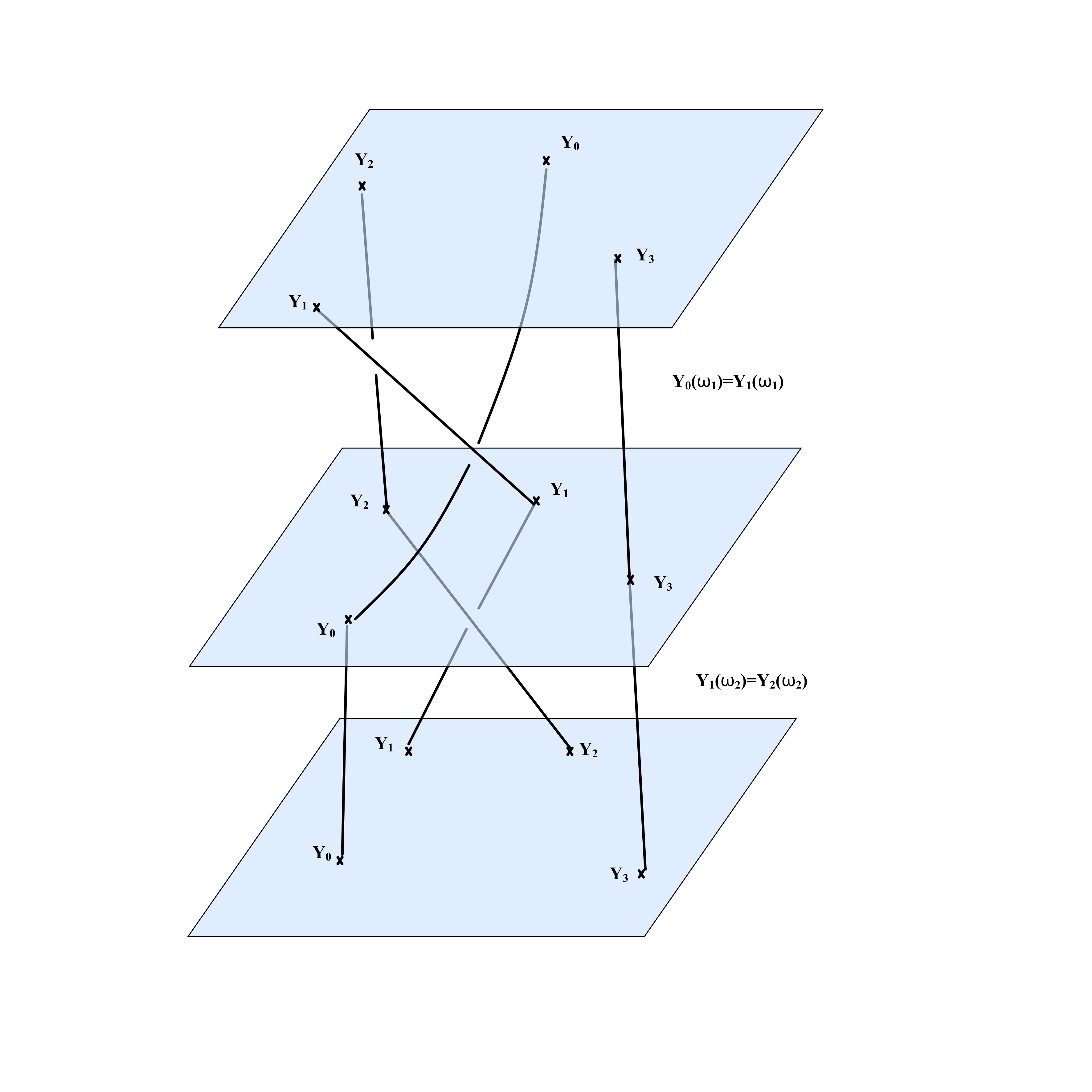}
\caption{ Position interchange of  the branches. }
\label{fig: coincidence}
\end{figure}

Figure \ref{fig: coincidence} illustrates  lemma \ref{Ytinterchange} for $n=4, p=1, r=5.$
In this case, the coincidence of roots $Y_0(\omega_1) = Y_1(\omega_1),$ $Y_1(\omega_2) = Y_2(\omega_2)$ happens.

\begin{figure}[H]
\centering
\includegraphics[scale=0.5]{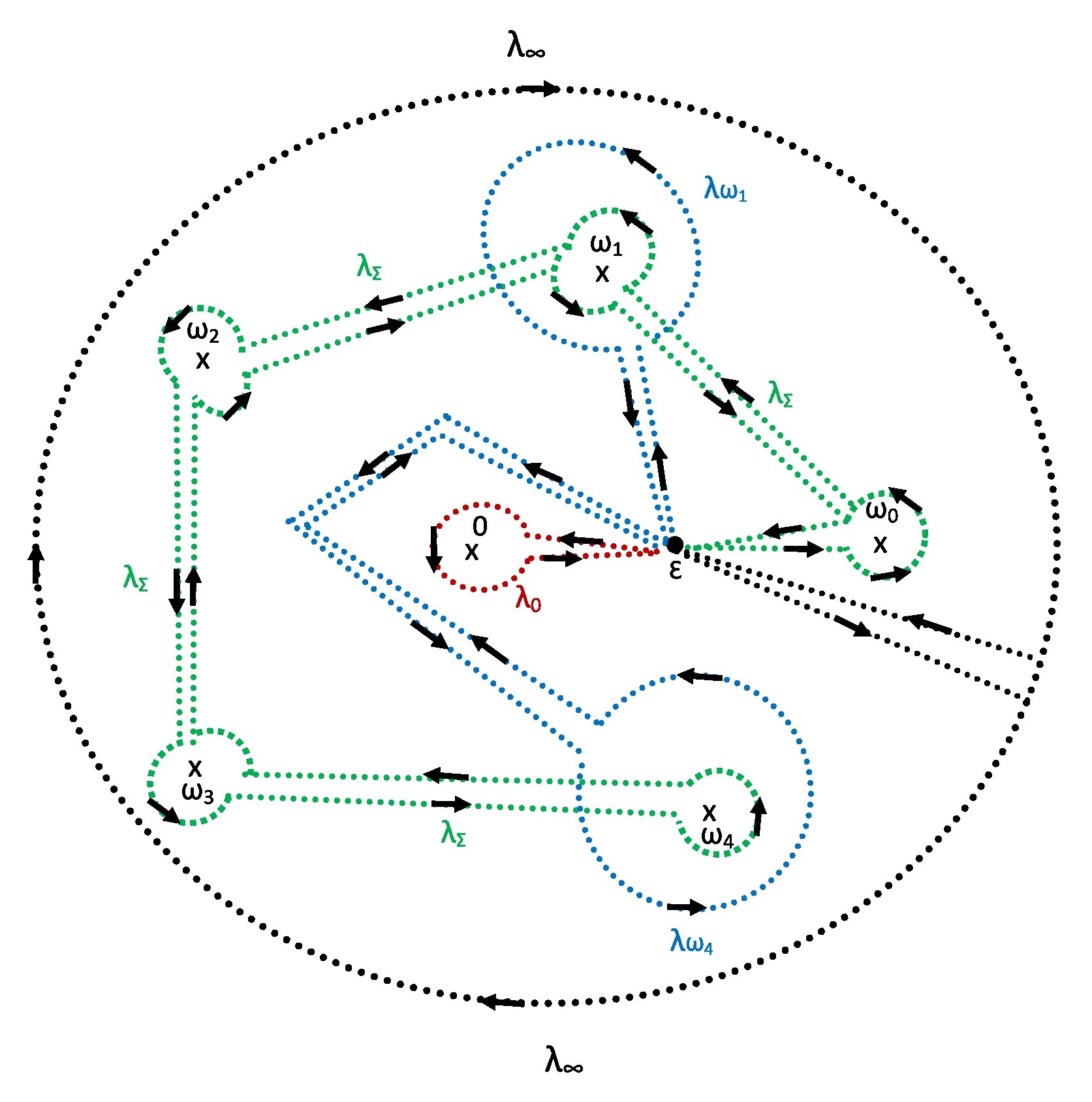}
\caption{ Path $\lambda_\Sigma$ }
\label{fig: LambdaSigma}
\end{figure}

Proposition  \ref{Ytcoincidence} and lemma \ref{Ytinterchange} furnish us with necessary data to describe the braid depicted by
the roots  $( Y_t(X))_{t=0}^{n-1}$  as $X$ makes a move along the path $\lambda_\Sigma.$

\begin{definition}\label{strandstwist} The $\alpha-$ twist with respect to a fixed point $Y=\Upsilon_0$ of $|{\mathcal K}|-$straight strands for 
${\mathcal K} \subset \{0, \cdots, n-1\},$
$$\{(s; Y_t(X) e(0 s))_{t\in {\mathcal K}} =(s; Y_t(X))_{t \in {\mathcal K}}; s \in [0,1]\},\ $$  is realized by the angular turn of $e(\alpha)$ as follows
\begin{equation} \label{alphatwist}
\{(s; (Y_t(X) - \Upsilon_0)e(\alpha s) + \Upsilon_0)_{t\in {\mathcal K}}; s \in [0,1]\}.
\end{equation} 
We shall denote this $\alpha-$ twist by ${\mathfrak R}^\alpha_{\Upsilon_0}[\mathcal K ].$ 
This notion is defined up to braid isotopy in the sense of Remark \ref{isotopy}.
\end{definition}

\begin{figure}[H]
\centering
\includegraphics[scale=0.3]{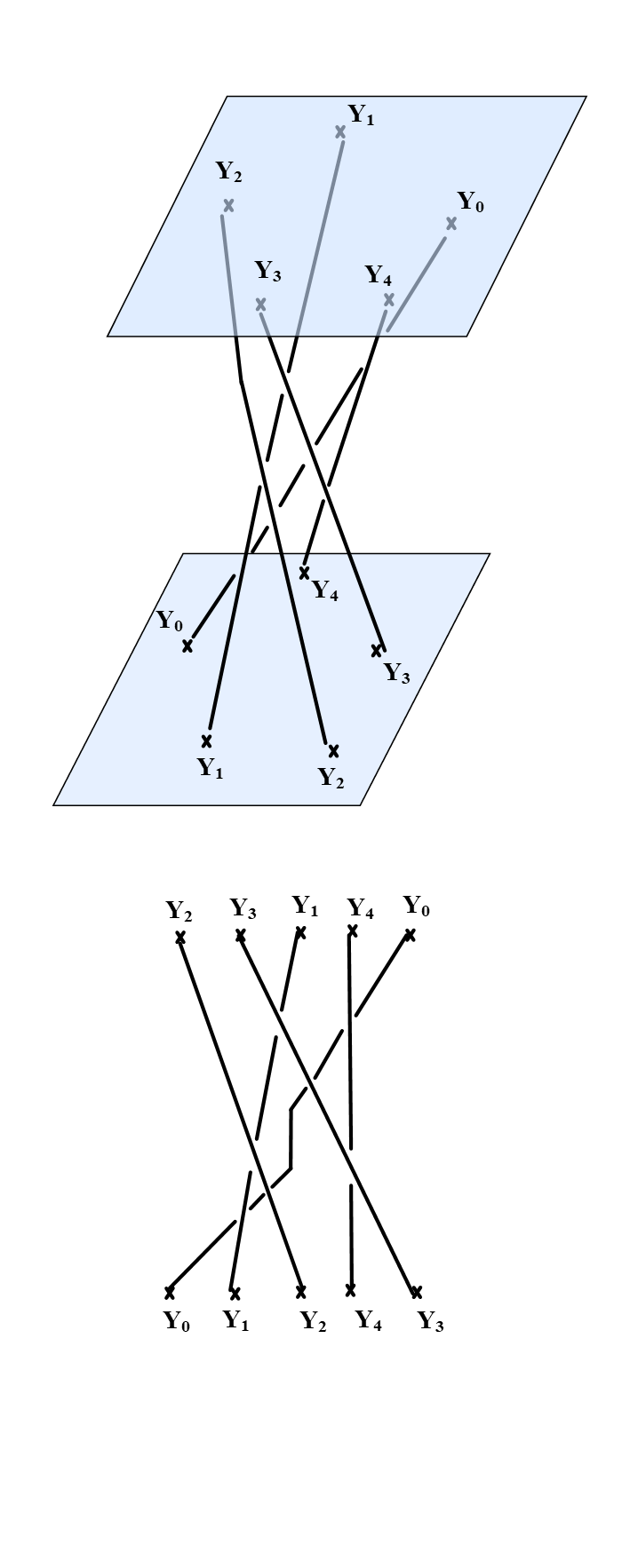}
\caption{ Twist ${\mathfrak R}^{\frac{2}{5}}_{\beta_0}[[0; 4] ]$ around  the barycenter $\beta_0$ of pentagon vertices: twice rotation of a pentagon \cite{MM}.
The lower part of the figure represents $\varpi({\mathfrak R}^{\frac{2}{5}}_{\beta_0}[[0; 4] ]) \in B_5$
after \eqref{lambdavarpi}, the corresponding Artin braid.  }
\label{fig: alphatwist}
\end{figure}

\begin{figure}[ht!]
\centering
\includegraphics[scale=0.3]{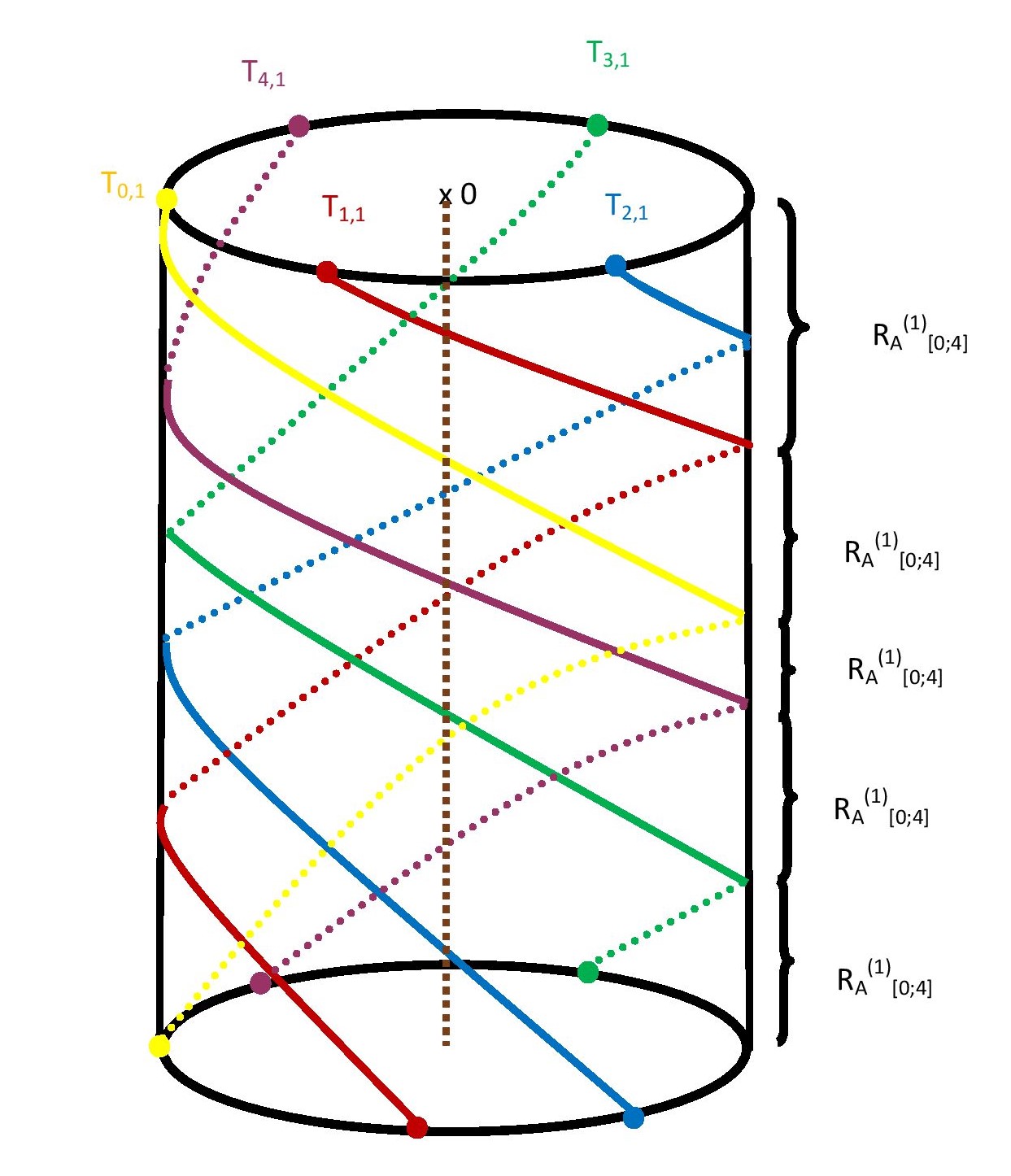}
\caption{ Twist ${\mathfrak R}^{1}_{0}[[0; 4] ]$ for $n=5:$ 5-rotation of a pentagon \cite{MM}. }
\label{fig: alphafulltwist}
\end{figure}

In the literature we find the notion of "rotating of punctures" \cite[Definition 1.4]{MM} that corresponds to
$\frac{1}{n}-$ twist for $n$ straight strands. The rotation of punctures is considered only for the case when the vertices of a polygon
make simultaneous moves in an anticlockwise direction so that the polygon itself remains the same. In our case, vertices may
make a movement of any angle. Thus the polygon itself may make a turn. See figures \ref{fig: alphatwist}, \ref{fig: alphafulltwist}.

We introduce the following notation well adapted to describe different kinds of $\alpha-$ twists.
\begin{definition}
For the standard Artin generators $\sigma_j, j \in [1;n-1]$ of $B_n,$ \cite{ArtinZopfe}, \cite{Birman}, we define
\begin{equation} \label{lambda}
\Lambda(\ell) = \sigma_\ell \sigma_{\ell-1} \cdots \sigma_1, \; \ell \in[0;n-1],
\end{equation} 
\begin{equation} \label{lambdabar}
\bar \Lambda(\ell) = \sigma_{n-\ell}  \cdots \sigma_{n-1}, \; \ell \in[1;n],
\end{equation} 
\begin{equation} \label{lambdaplus}
\Lambda^+(k) = \sigma_{1} \cdots \sigma_{k}, \; k \in[0;n-1].
\end{equation} 

\begin{equation} \label{lambdaminus}
\Lambda^-(j) =\sigma_{n-1} \cdots \sigma_{n-j}, \; j \in[1;n].
\end{equation} 
\end{definition}
 Here we use the convention
 $$ \Lambda(0) =\bar \Lambda(n) = \Lambda^+(0) =  \Lambda^-(n) = id.$$
 
The following proposition gives a translation from the language of angular twists to
that of Artin braids. It shall be, however, noted that this correspondence does not
necessarily imply an homomorphism from twists to Artin braids. That is to say, if $\alpha_j-$twists correspond to
Artin braids $\Lambda_j, j =1,2,$  $(\alpha_1+\alpha_2)-$twist does not necessarily corresponds to
$\Lambda_1 \cdot \Lambda_2.$ This can be seen from an
inequality like 
$$\left\lfloor \frac{3}{5} \right\rfloor + \left\lfloor \frac{2}{5} \right\rfloor \not = \frac{3}{5} + \frac{2}{5}.$$ 
This correspondence can be understood as a projection of  movement of strands in three dimensional space onto a two dimensional
plane. Thus it depends on the choice of the projection direction.  We can deduce
Artin braid description of algebraic functions from their twist description, but the converse is not true in general. 
In the sequel, we shall denote 
\begin{equation} \label{varpitau} 
 \varpi({\mathfrak R}^\alpha_{\Upsilon_0}[\mathcal K ]) \in B_n
\end{equation}
the projection of  a twist ${\mathfrak R}^\alpha_{\Upsilon_0}[\mathcal K ], {\mathcal K} \subset [0; n-1],$ onto an Artin braid.
For each $\varpi,$ we choose a fixed projection direction.

\begin{proposition} \label{twistbraid}
For $\alpha \in [\frac{k}{n}, \frac{k+1}{n}), k \in \Z_{\geq 0},$ we denote by $k^-$ a non-negative integer equal to or smaller than
$k-1.$  
An Artin braid corresponding to the twist ${\mathfrak R}^\alpha_{\beta_0}[[0; n-1]]$ around the barycenter $\beta_0$ of a $n-$polygon
admits the following expression,
\begin{equation} \label{alphalambda} 
\varpi({\mathfrak R}^\alpha_{\beta_0}[[0; n-1]]) = \Lambda(\ell_1) \cdots \Lambda (\ell_{u}) \left( \Lambda(n-1) \right)^{k^-} \Lambda^-(\ell^-_{u^-}) \cdots \Lambda^-(\ell^-_{1}) .
\end{equation}
by \eqref{lambda},  \eqref{lambdaminus} with
$0 \leq \ell_1 \leq  \cdots \leq \ell_{u} < n-1,$
$1 \leq \ell_1^- \leq  \cdots \leq \ell_{u^-}^- \leq n.$ In 
\eqref{alphalambda} no factor $\Lambda^-(n-1)$ appears.

In a similar manner, we have also
\begin{equation} \label{alphalambda1} 
\varpi({\mathfrak R}^\alpha_{\beta_0}[[0; n-1]]) = {\bar \Lambda}(\bar \ell_{1}) \cdots {\bar \Lambda}(\bar \ell_{\bar v})
   (\Lambda^+(n-1))^{k^-} \Lambda^+(\ell_{v^+}^+) \cdots \Lambda^+ (\ell_1^+) .
\end{equation}
in terms of   \eqref{lambdabar}, \eqref{lambdaplus} with
$0 \leq \ell_1^+ \leq  \cdots \leq \ell_{v^+}^+ < n-1,$
$1 \leq \bar \ell_1\leq  \cdots \leq \bar \ell_{\bar v} \leq n.$
 In 
\eqref{alphalambda1} no factor ${\bar \Lambda}(n-1)$ appears.
We notice that this decomposition expression is not unique. 
\end{proposition}

See Figures \ref{ArtindiagramQuinticeqn1}, \ref{ArtindiagramQuinticeqn2} for a concrete example of
the above Proposition.  

The relation \eqref{varpitau} mentioned above can be reformulated into the following diagram. Here
$\lambda_\ast \in {\mathcal T}_n$ denotes a groupoid  of twists that are induced from a loop 
$ \lambda \in \pi_1 (\C^\ast \setminus \Sigma, \epsilon).$ This is an example of semi-direct product groupoids
treated in \cite[Example 3]{Brown}. The groupoid of twists ${\mathcal T}_n$ can be obtained as a composition of elementary $\alpha-$twists of Definition \ref{strandstwist} defined up to braid isotopy in the sense of Remark \ref{isotopy}.
 Two twists ${\mathfrak R}_j, j = 1,2, $ with fixed ends ${\mathfrak R}_j\cap ( \{0,1\} \times \C),$
are considered to be equivalent in ${\mathcal T}_n$ 
if two open sets $ \bigl((0,1) \times \C \bigr) \setminus {\mathfrak R}_j$ that are complement to $n-$strands  in $\bigl( (0,1) \times \C \bigr)$ are homeomorphic one to another.  
\begin{equation}
 \begin{array}{ccccc}
              \pi_1 (\C^\ast \setminus \Sigma, \epsilon) & \longrightarrow &{\mathcal T}_n & \overset{\varpi}{\longrightarrow} & B_n\\
              \rotatebox{90}{$\in$}& & \rotatebox{90}{$\in$}& &\rotatebox{90}{$\in$}\\ 
\lambda&\longmapsto & \lambda_\ast  & \longmapsto &  \varpi (\lambda_\ast)
             \end{array} .
\label{lambdavarpi}
\end{equation}
The first arrow gives a groupoid homomorphism between the fundamental group and ${\mathcal T}_n$, while the second arrow $\varpi$
is a projection. 
\begin{example}\label{Garside}
The Garside element $\Delta^2$ that generates the center $Z(B_n) \subset B_n$
is realized by $\alpha= 1$ twist or $n-$rotation of punctures after \cite[Definition 1.4]{MM}. 
That is to say, all strands make a $2 \pi$ turn in \eqref{alphatwist}.

The $1/2$ twist $\Delta$ admits a representation expressed by  \eqref{lambda}, \eqref{lambdaplus}
$$ \Delta = \Lambda(1) \cdots \Lambda(n-1) =  \Lambda^+(n-1) \cdots \Lambda^+(1).$$ 
\begin{figure}[H]
\centering
\includegraphics[scale=0.5]{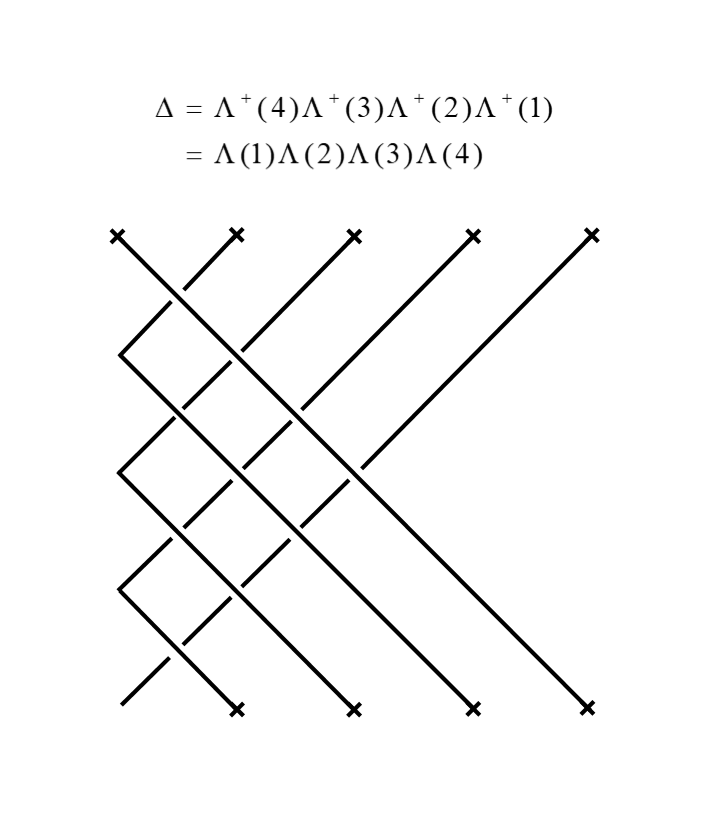}
\caption{ $\frac{1}{2}-$ twist $\Delta$ }
\label{fig: Delta}
\end{figure}
\end{example}

Further we denote by  $\lambda_\ast$ the monodromy action on $Y_t(X), t \in [0;n-1]$ provoked by a loop $\lambda$ lying in the open set $\C \setminus (\Sigma \cup \{0\}).$

We shall study monodromy actions induced by the following loops depicted in the Figure \ref{fig: LambdaSigma}.
For all of them we choose a starting point $\epsilon \in \{X; |\frac{X^N}{R}|<1\}.$ The loop makes an anticlockwise turn only around the indicated branching points and no other branching points. 

\vspace{0.5pc}
$\bullet$  $\lambda_0:$ a loop starting from $\epsilon$ that makes an anticlockwise turn around $X=0.$ 

\vspace{0.5pc}
$\bullet$  $\lambda_\Sigma:$ a loop starting from $\epsilon$ that makes an anticlockwise turn around $\Sigma$ \eqref{branching points}. 

\vspace{0.5pc}
$\bullet$  $\lambda_\infty:$ a loop starting from $\epsilon$ that makes an anticlockwise turn around $X=\infty.$ 

\vspace{0.5pc}
$\bullet$  $\lambda_{\omega_\ell}:$ a loop starting from $\epsilon$ that makes an anticlockwise turn around $X=\omega_\ell \in \Sigma, \ell \in [0;N-1],$  \eqref{branching points}.

\begin{theorem}\label{AlgBraid}
We establish the following
braid group representation of the monodromy of algebraic functions $Y_t(X), t \in [0; n-1]$ as $X$ makes movement along a loop starting from a   point $\epsilon \in \{X; |\frac{X^N}{R}|<1\}$ that avoids branching points $\Sigma \cup \{0, \infty\}. $
Here we use the notation \eqref{alphatwist} to represent an angular twist. 

(I) The monodromy $(\lambda_0)_\ast$ is given by a $\frac{r-g}{p}-$ twist ${\mathfrak R}^{\frac{r-g}{p}}_{0}[[0; p-1]]$  (resp. a $\frac{g}{q}-$ twist $ {\mathfrak R}^{\frac{g}{q}}_{0}[[p; n-1]]$) 
for $p-$solutions \eqref{psolns} (resp. q- solutions\eqref{qsolns}).

(II)
The monodromy $(\lambda_\Sigma )_\ast$ is given by
a $-\frac{N}{np}$ twist ${\mathfrak R}^{-\frac{N}{np}}_{0}[[p; n-1]]$  (resp. a $\frac{N}{nq}$ twist ${\mathfrak R}^{\frac{N}{nq}}_{0}[[p; n-1]]$ )  for $Y_{u}(X), u \in [0; p-1].$ (resp.  $Y_{v}(X), v \in [p; n-1]$).
  The local braid resulting from  $\lambda_{\omega_{\bar \ell}}$ with $\omega_{\bar \ell} \in \Sigma$ can be described after the rule \eqref{YtYt+1coincidence} of Lemma \ref{Ytinterchange} i.e.
 $Y_{t}$ and $Y_{t'}$  achieve a $\frac{1}{2}-$twist   ${\mathfrak R}^{\frac{1}{2}}_{c_{\bar \ell}}[\{t, t'\}]$   for $c_{\bar \ell}=Y_t(\omega_{\bar \ell})=Y_{t'}(\omega_{\bar \ell})$ where the index $t$ is determined by
\begin{equation}\label{lprt1} \bar\ell \equiv p(rt+1) \quad \mod n,
    \end{equation}
and the index $t'$ is determined by
\begin{equation}\label{rtt}r(t'-t) \equiv 1 \mod n.
\end{equation}
The monodromy $(\lambda_\Sigma)_\ast$ is braid isotopic to the composition $\prod_{\ell=0}^{N-1}(\lambda_{\omega_\ell})_\ast$
in the sense of Remark \ref{isotopy}.

(III) The monodromy $(\lambda_\infty)_\ast$ is given by a $(-\frac{r}{n})-$ twist 
${\mathfrak R}^{-\frac{r}{n}}_{0}[[0; n-1]]$ 
  for all solutions \eqref{psolns}, \eqref{qsolns}.
\end{theorem}

\begin{proof}
(I) The  series expansions for  $p-$solutions \eqref{psolns} (resp. $q-$solutions \eqref{qsolns}) show that the monodromy $(\lambda_0)_\ast$ 
induces $Y_u(X) \rightarrow e(\frac{r-g}{p}) Y_u(X), u \in [0;p-1]$  
(resp. $Y_{v}(X) \rightarrow e(\frac{g}{q}) Y_{v}(X), v \in [p;n-1]$).

(III) The asymptotic expansion \eqref{Ytinfty} valid near $ X =\infty,$ shows that 
the monodromy action along a path homotopic to $\lambda_\infty$ that turns around $X=\infty$ in an anticlockwise way yields
$- (\frac{r}{n})-$ twist with respect to $Y=0$ of $(Y_t(X))_{t=0}^{n-1}.$ 

(II)
The statement (III) above indicates that the composition of monodromy actions $(\lambda_0 \lambda_\Sigma )_\ast$  gives rise to a $\frac{r}{n} -$ twist with respect to $Y=0$ for both of $p-$solutions \eqref{psolns} and $q-$solutions \eqref{qsolns}. 
The monodromy $(\lambda_\Sigma )_\ast$ can be extracted from this $\frac{r}{n} -$ twist and the monodromy
$(\lambda_0)_\ast$ described in (I).
In other words, the monodromy $(\lambda_\Sigma)_\ast$ corresponds to a
$( \frac{r}{n} - \frac{r-g}{p})  -$ twist with respect to $Y=0$ for $p-$solutions \eqref{psolns}
(resp.  $(\frac{r}{n}- \frac{g}{q})-$ twist with respect  to $Y=0$ for $q-$solutions \eqref{qsolns}).

Proposition \ref{Ytcoincidence}, Lemma \ref{Ytinterchange} give a description of each local twist provoked by $\lambda_{\omega_{\bar \ell}}$
that is achieved by $( Y_t(X))_{t=0}^{n-1}$ at $\omega_{\bar \ell}, \bar \ell \in [0, N-1].$
That is to say, $(\lambda_\Sigma)_\ast$ is braid isotopic to the composition $\prod_{\ell=0}^{N-1}(\lambda_{\omega_\ell})_\ast$
in the sense of Remark \ref{isotopy}. 
\end{proof}

The Lemma 3.5 of  \cite[3.2 Tropical construction of coarse braids]
{EsterovLangBraid} implies our Theorem \ref{AlgBraid}, (III)  above.  

One may wonder why a negative twist around the origin $Y=0$ appears as a result of
the monodromy $(\lambda_\Sigma)_\ast$ for $q-$solutions
in Theorem \ref{AlgBraid}, (II). We explain the situation with the example of a quintic equation.

\begin{example}\label{n3p1r4}
To illlustrate Theorem \ref{AlgBraid}, let us look at an example of a quintic equation
$$ f(X,Y) = Y^5- X^2 Y^3 + X^7,$$
with $N=4, \Sigma =\{ \omega_\ell = e(\frac{\ell}{4})(\frac{2^2 3^3}{5^5})^{\frac{1}{4}} \}_{\ell =0}^3.$
We have three  $(p=3)-$solutions  $Y_0(X), Y_1(X), Y_2(X)$  and two $(q=2)-$solutions $Y_3(X), Y_4(X)$.

The monodromy $(\lambda_0)_\ast$ yields $\frac{5}{3}-$twist ${\mathfrak R}^{\frac{5}{3}}_{0}[[0; 2]]$  for 
 $(p=3)-$solutions (resp. $1-$ twist ${\mathfrak R}^{1}_{0}[[3; 4]]$ for $(q=2)-$solutions).

The interchanging roots provoked by the monodromy $(\lambda_\Sigma)_\ast$
are given by the following list:
$$  Y_2(\omega_0) = Y_0(\omega_0),  Y_1(\omega_1) = Y_3(\omega_1),  Y_2(\omega_2) = Y_4(\omega_2),  
Y_3(\omega_3) = Y_0(\omega_3).   $$ 
Even though each of two solutions make an anticlockwise turn with respect to one another, along $\lambda_\Sigma$
$(p=3)-$solutions achieve an argument change by $-\frac{8}{15}\pi$ (i.e. a twist ${\mathfrak R}^{-\frac{4}{15}}_0[[0;2]]$.)
while  the argument change of the $(q=2)-$solutions  is $\frac{4}{5}\pi.$
In total, the monodromy $(\lambda_0 \cdot \lambda_\Sigma)_\ast$
yields $(\frac{7}{5})-$ twist ${\mathfrak R}^{\frac{7}{5}}_0[[0;4]]$ for all solutions. In this way the composition, of 
$(\lambda_0 \cdot \lambda_\Sigma)_\ast$ and $(\lambda_\infty)_\ast$
gives rise to a trivial twist.
\begin{figure}[H]
\centering
\includegraphics[ totalheight=18cm]{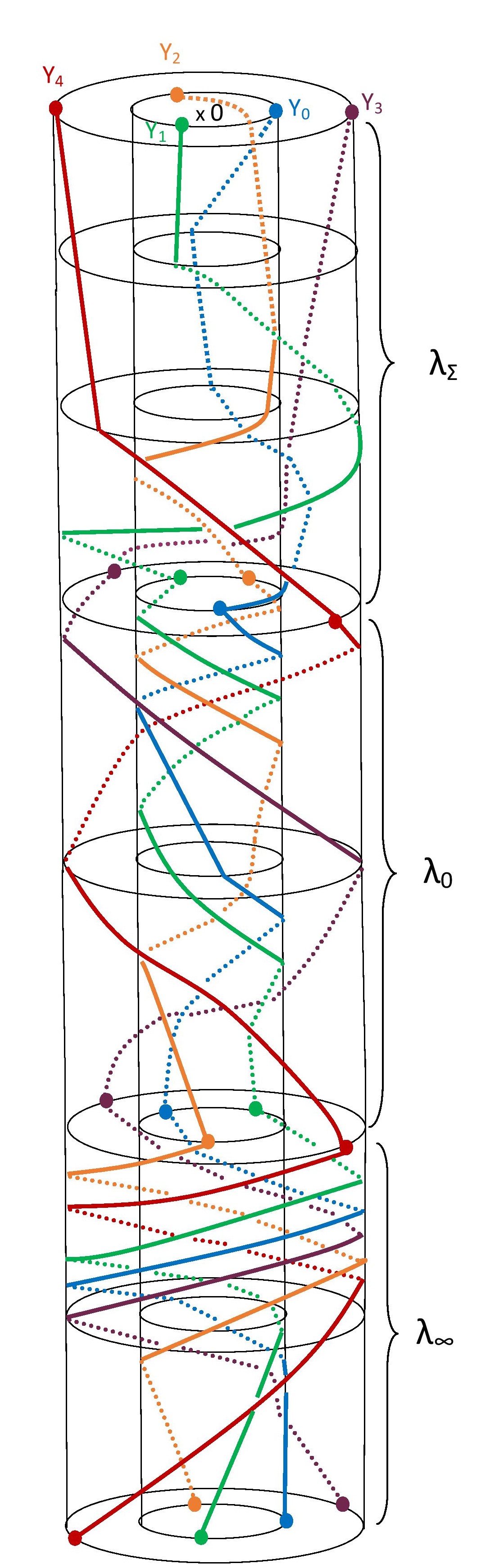}
\caption{ Global braid for $Y^5- X^2 Y^3 + X^7.$  }
\label{fig: quinticpq}
\end{figure}

\begin{figure}[H]\label{ArtindiagramQuinticeqn1}
    \centering
     \begin{tikzpicture}[thick,scale=0.6, every node/.style={scale=0.6}]
\pic[ 
  rotate=0,
  braid/.cd,
  every strand/.style={ultra thick},
  strand 1/.style={red},
  strand 2/.style={purple},
  strand 3/.style={orange},
  strand 4/.style={blue},
  strand 5/.style={green},
] {braid= s_1 s_3 s_4 s_2 s_3 s_1 s_4 s_2 s_3 s_4 s_1 s_2 s_3 s_4 s_1 s_2 s_3 s_4 s_1 s_2 s_3 s_4 s_1 s_2 s_3 s_1 } ;
\node (1) at (0,0.2) {$\color{red}Y_1$};
\node (2) at (1,0.2) {$\color{purple}Y_2$};
\node (3) at (2,0.2) {$\color{orange}Y_4$};
\node (4) at (3,0.2) {$\color{blue}Y_0$};
\node (5) at (4,0.2) {$\color{green}Y_3$};
\node (6) at (-1,-6.5) {};
\node (7) at (6,-6.5) {};
\draw[dashed](6)--(7);
 \node (8) at (7,-3.25){}; 
 \node (9) at (10,-15.25){}; 
\end{tikzpicture}
\caption{Artin braid diagram $\varpi\bigl((\lambda_{\Sigma})_{\ast}\bigr)\varpi\bigl((\lambda_{0})_{\ast}\bigr)$ $= \Lambda(1) \Lambda(3) (\Lambda(4))^5 \Lambda^-(2)$ $=\bar \Lambda(2)  (\Lambda^+(4))^5 \Lambda^+(3) \Lambda^+(1)$ for $ Y^5-X^2Y^3+X^7=0.$  }
    
\end{figure}
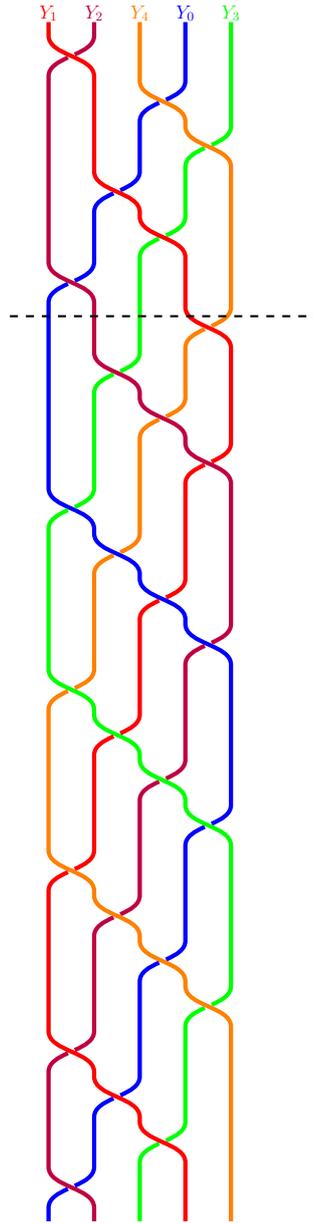


\begin{figure}[H]
    \centering
     
     \begin{tikzpicture}[thick,scale=0.6, every node/.style={scale=0.6}]
\pic[ 
  rotate=0,
  braid/.cd,
  every strand/.style={ultra thick},
  strand 1/.style={blue},
  strand 2/.style={purple},
  strand 3/.style={green},
  strand 4/.style={red},
  strand 5/.style={orange},
] {braid= s_1^{-1} s_3^{-1} s_2^{-1} s_1^{-1} s_4^{-1} s_3^{-1} s_2^{-1} s_1^{-1} s_4^{-1} s_3^{-1} s_2^{-1} s_1^{-1} s_4^{-1} s_3^{-1} s_2^{-1} s_1^{-1} s_4^{-1} s_3^{-1} s_2^{-1} s_1^{-1} s_4^{-1} s_3^{-1} s_2^{-1} s_1^{-1} s_4^{-1} s_3^{-1}  } ;
\node (1) at (0,0.2) {$\color{blue}Y_0$};
\node (2) at (1,0.2) {$\color{purple}Y_2$};
\node (3) at (2,0.2) {$\color{green}Y_3$};
\node (4) at (3,0.2) {$\color{red}Y_1$};
\node (5) at (4,0.2) {$\color{orange}Y_4$};
\node (6) at (-1,-6.5) {};
\node (7) at (6,-6.5) {};
 \node (8) at (14,-13.25) {
 }; 
\end{tikzpicture}

    \caption{Artin braid diagram $\varpi \bigl((\lambda_{\infty})_{\ast}\bigr)$  $= \Lambda^-(2)^{-1}(\Lambda(4))^{-5} \Lambda(3)^{-1} \Lambda(2)^{-1}  $ $= \Lambda^+(1)^{-1} \Lambda^+(3)^{-1}   (\Lambda^+(4))^{-5}\bar \Lambda(2)^{-1} $ for $ Y^5-X^2Y^3+X^7=0.$ }
    \label{ArtindiagramQuinticeqn2}
\end{figure}
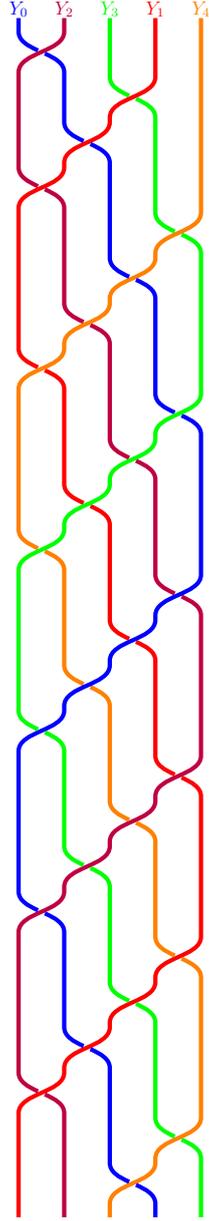
\end{example}

As a corollary to our Theorem \ref{AlgBraid}, we establish Theorem 1 of 
\cite{EsterovLangBraid} and \cite[Corollary 3.6]{KatoNoumi}. Even though we have far less loops on the moduli space of algebraic equations at our disposal than in \cite{EsterovLangBraid} where loops are chosen in $\C^3$, we conclude that the induced braid monodromy group of roots  $B_n(f) = \varpi ({\mathcal T}_n),$ \eqref{exactseq} does recover the full braid group under \eqref{gcdnN}
that is a condition inessential in \cite[Theorem 1]{EsterovLangBraid}.

\begin{corollary}\label{BnfBn}
Under the conditions $h.c.f. (n,p) = h.c.f. (n,N) =1,$ \eqref{gcdnN}
we have the following isomorphisms for the braid monodromy $B_n(f)$ (resp. Galois group $Gal(f)$) of $\eqref{f(X,Y)0}$
\begin{equation}\label{Bnf}
B_n(f) \simeq B_n, \;\; Gal(f) \simeq {\mathfrak S}_n.
\end{equation}

\end{corollary}
\begin{proof}
First of all, we recall the definition of the algebraic functions \eqref{Capital} that would indicate a relation  $(\lambda_\infty)_\ast^{-j}Y_{t+j}(X) = Y_t(X). $
Among all roots subject to the monodromy action $(\lambda_\infty)_\ast^{-j}$, the pair of roots
$Y_{t+j}, Y_{t'+j}$ interchange their relative position at $\omega_{\bar \ell}$
for the indices $(t,t') \in [0;n-1]^2, \bar \ell \in [0;N-1]$
determined by the formulae \eqref{lprt1}, \eqref{rtt}.
This means that the monodromy action 
$(\lambda_\infty^{-j} \lambda_{\omega_{\bar \ell}} \lambda_\infty^{j})_\ast$
produces the interchange of roots  $Y_{t+j}, Y_{t'+j}$ and leaves $Y_v, v \in [0; n-1 ] \setminus \{t+j, t'+j\}$ invariant. 

A theorem that goes back to E.Artin \cite[Satz 1]{ArtinZopfe}, \cite[Theorem 1.8]{Birman} tells us that all interchanges of pair of roots with indices $\{(t+j, t'+j)\}_{j=0}^{n-1}$ generate the braid group $B_n.$

We deduce the statement on the Galois group $Gal(f)$ of the roots of $\eqref{f(X,Y)0}$ from the exact sequence \eqref{exactseq}.
\end{proof}

{
\center{\section{Non-coprime $(n,p)$ case}\label{noncoprime}}
}
Now we consider the equation
\begin{equation}\label{f(X,Y)1}
f(X,Y)=Y^{mn}-X^{g}Y^{mp}+X^{r}=0.
\end{equation}
with $m>1,$ $(n,p)$ coprime.
In this situation, the description of the monodromy representations $(\lambda_0)_\ast, (\lambda_\infty)_\ast$ 
are parallel to the case \eqref{f(X,Y)0} for $m=1$ in \eqref{f(X,Y)1}. A new feature appears in the monodromy 
$(\lambda_{\Sigma} )_\ast.$

To describe the braid monodromy of the solutions to \eqref{f(X,Y)1},
we regard it as a $m-$covering of the following master equation
\begin{equation}\label{F(X,Y)1}
F(X,T)=T^{n}-X^{g}T^{p}+X^{r}=0,
\end{equation}
with solutions $\{T_t(X)\}_{t=0}^{n-1}.$
Then we have the solutions to \eqref{f(X,Y)1} in the form
\begin{equation}\label{Ttm}
 Y_t^{(j)}(X) = e(\frac{j}{m})  T_t^{1/m}(X)\;\; t\in [0;n-1],  j\in [0;m-1].
\end{equation}

The branching points of the solutions to \eqref{f(X,Y)1}
are given by $\Sigma$ \eqref{branching points}.

In applying the result of Theorem\ref{AlgBraid} to \eqref{F(X,Y)1}, we get immediately the description of $(\lambda_0)_\ast,$  $(\lambda_{\Sigma} )_\ast$ and $(\lambda_\infty)_\ast.$

The only novelty appears in the braid monodromy around one of
branching points $\Sigma$ \eqref{branching points}.

In order to describe the braid monodromy at a branching point $\omega_\ell \in \Sigma,$
$\ell \in [0; N-1],$
we build $n-$tuple of curved polyhedra called "roof tiles."
For a fixed $t \in [0;n-1],$ we consider a union of $(m-1)-$edges of a polygon with vertices 
\begin{equation}\label{Tjell}
 \{  Y^{(j)}_t(\tau_\ell (\theta))\}_{j=0}^{m-1},
 \end{equation}
defined over a continuous closed curve $\tau_\ell :[0, 1] \ni \theta \longrightarrow \C \setminus (\Sigma \cup \{0\}).$
 Especially we are interested in a loop starting from
$\tau_\ell(0) = \tau_\ell(1)= \epsilon$ that makes an anticlockwise turn around $\omega_\ell \in \Sigma,$ $\ell \in [0; N-1].$

We put the above mentioned union of edges, i.e. a polygon lacking one side, into a domain $\Delta_t(\theta), t \in [0; n-1],$
that is homeomorphic to a disc. The vertices \eqref{Tjell} belong to $\Delta_t(\theta)$ for each fixed $t.$ 
We require that $\Delta_t(0) \cap \Delta_{t'}(0) = \emptyset , \forall t \not = t'.$
Each of $\Delta_t(\theta)$ is obtained as a result of rotation of $\Delta_t(0)$ 
that is achieved by the movement of  $\theta \in [0,1].$ Thus we construct a $n-$tuple of  curved polyhedra that we call ''roof tiles''
\begin{equation}\label{Deltabar}
\bar \Delta_t = \bigcup_{\theta \in [0,1]}\left(\theta, \cup_{t=0}^{n-1} \Delta_t(\theta) \right)  \subset [0,1] \times \C,
 \end{equation}
with $t \in [0; n-1].$

\begin{figure}[H]
\centering
\includegraphics[scale=0.6]{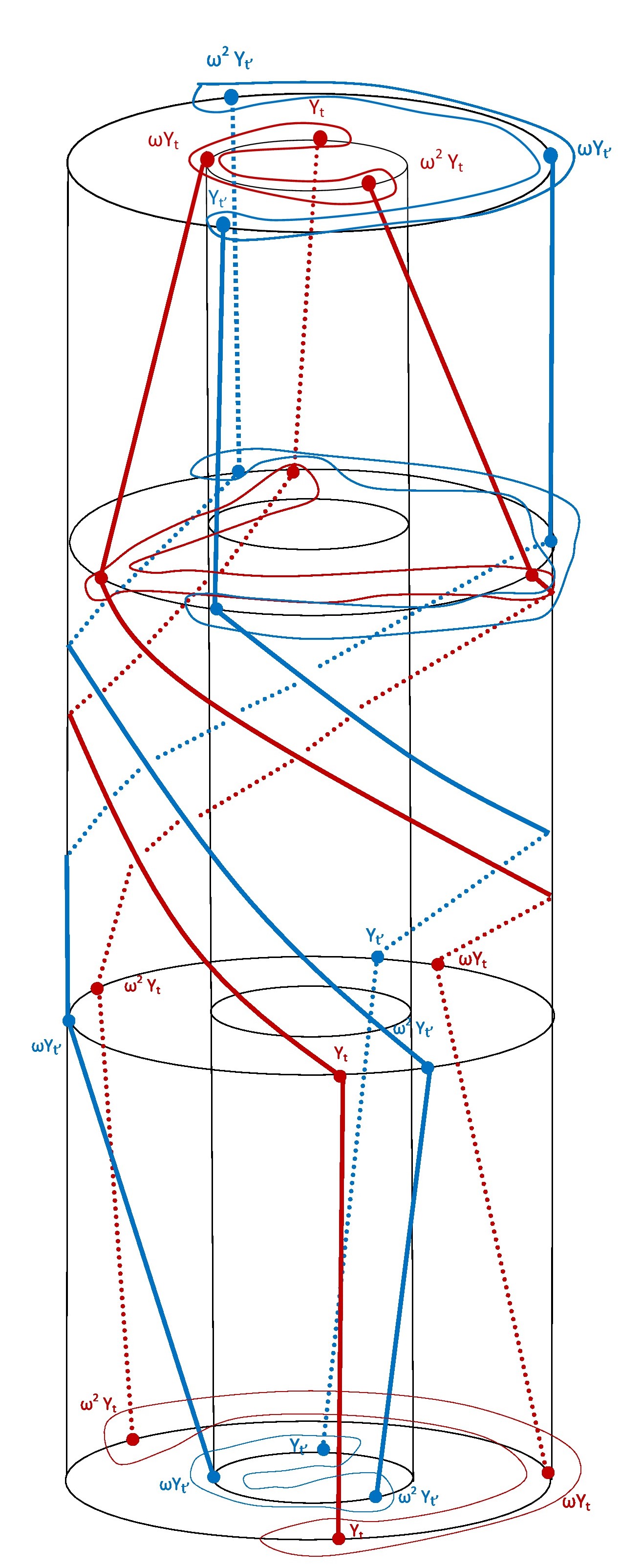}
\caption{Roof tiles $\bar \Delta_0,$  $\bar \Delta_1$ for $m=3$. }
\label{fig:Tuiles}
\end{figure}

Theorem \ref{AlgBraid}, together with the definition of the roots
$Y_t^{(j)}(X), t\in [0;n-1],  j\in [0;m-1],$ \eqref{Ttm}
yields the following.
\begin{theorem}\label{AlgBraidnoncoprime}
We establish the following
braid group representation of the monodromy of algebraic functions
 $Y_t^{(j)}(X), t \in [0; n-1], j \in [0; m-1]$ \eqref{Ttm} as $X$ makes movement along a loop starting from a   point $\epsilon \in \{X; |\frac{X^{N}}{R}|<1\}$ that avoids branching points $\Sigma \cup \{0, \infty\}. $

(I) The monodromy $(\lambda_0)_\ast$ is given by a $\frac{r-g}{mp}-$ twist ${\mathfrak R}^{\frac{r-g}{mp}}_{0}[[0; p-1] \times [0; m-1]]$ (resp. $\frac{g}{mq}-$ twist ${\mathfrak R}^{\frac{g}{mq}}_{0}[[p; n-1]\times [0; m-1]]$) for the $m-$th root \eqref{Ttm} of  $p-$solutions \eqref{psolns}   with $t \in [0;p-1]$ (resp. the $m-$th root of q- solutions\eqref{qsolns} with $t \in [p;n-1]$).

(II)
The monodromy $(\lambda_{\Sigma })_\ast$ is given by
a $(-\frac{N}{mnp})-$twist ${\mathfrak R}^{-\frac{N}{mnp}}_{0}[[0; p-1]\times [0; m-1]]$ (resp. $\frac{N}{mnq}-$ twist ${\mathfrak R}^{\frac{N}{mnq}}_{0}[[p; n-1]\times [0; m-1]]$) for $Y_{u}(X), u \in [0; p-1].$ (resp.  $Y_{v}(X), v \in [p; n-1]$).
 The local braid resulting from  $\lambda_{\omega_{\bar \ell}}$ with $\omega_{\bar \ell} \in \Sigma$ can be described after the rule \eqref{YtYt+1coincidence} applied to $ \Delta_{t}(\theta)$ and $ \Delta_{t'}(\theta)$ for $\theta = 0,1.$ 
The movement of $ \Delta_{t}(\theta)$ (resp. $ \Delta_{t'}(\theta)$), $\theta \in [0,1],$ is restricted inside the roof tile  $\bar \Delta_{t}$ (resp.  $\bar \Delta_{t'}$).
In other words, for each $j \in [0;m-1],$ $Y^{(j)}_{t}(X)$ and $Y^{(j)}_{t'}(X)$ achieve a $\frac{1}{2}-$ twist ${\mathfrak R}^{\frac{1}{2}}_{c_{\bar \ell} }[\{t, t'\}\times \{j\}]$ 
$$c_{\bar \ell}=  Y^{(j)}_{t}(\omega_{\bar \ell}) = Y^{(j)}_{t'}(\omega_{\bar \ell}), $$ 
the center of symmetry of roof tiles while the angular ordering
of  $Y^{(j)}_{t}(X)$ (resp. $Y^{(j)}_{t'}(X)$), $j \in [0;m-1],$ inside the roof tile $\bar \Delta_{t}$ (resp.  $\bar \Delta_{t'}$)
is preserved during the  movement $X=\tau_{\bar \ell}(\theta), \theta \in [0,1]. $
That is to say,
$$ Arg \left( Y^{(0)}_{t}(\tau_{\bar \ell}(\theta)) - c_{\bar \ell} \right) < Arg \left( Y^{(1)}_{t}(\tau_{\bar \ell}(\theta)) - c_{\bar \ell} \right)< \cdots < Arg \left( Y^{(m-1)}_{t}(\tau_{\bar \ell}(\theta)) - c_{\bar \ell}\right), $$
$$ Arg \left( Y^{(0)}_{t'}(\tau_{\bar \ell}(\theta)) - c_{\bar \ell} \right) < Arg \left( Y^{(1)}_{t'}(\tau_{\bar \ell}(\theta)) - c_{\bar \ell} \right)< \cdots < Arg \left( Y^{(m-1)}_{t'}(\tau_{\bar \ell}(\theta)) - c_{\bar \ell}\right), $$
$\forall \theta\in[0,1].$

(III) The monodromy $(\lambda_\infty)_\ast$ is given by a $(-\frac{r}{mn})-$ twist
${\mathfrak R}^{-\frac{r}{mn}}_{0}[[0; n-1]\times [0; m-1]]$
  for all solutions \eqref{Ttm}.
\end{theorem}

We remark here that the movement of  $ \Delta_{t}(\theta)$ and $ \Delta_{t'}(\theta),$ $\theta \in [0;1]$
does not represent an example of "convex twist" introduced in \cite[Defintion 3.1]{MM}.
In the do-si-do dance similar to the convext twist, two dancers approach each other and circle back to back. The direction they face is unchanged.
In our case, the two dancers face one to another all the time and never take the position of back to back.

As a corollary to our Theorem \ref{AlgBraidnoncoprime}, we establish Corollary 2 of 
\cite{EsterovLangBraid} on the wreath product of the Artin braid group if $h.c.f. (n,N) =1.$ 
This means that our result recovers \cite[Theorem 1.5]{EsterovLangWreath} on the wreath product structure of the Galois group for \eqref{f(X,Y)1}. In other words,  the Galois
group $Gal(f)$ for \eqref{f(X,Y)1} has order $m^{n-1} n!$ as it has been announced in \cite[Corollary 4.6]{KatoNoumi}.

Let us recall the definition of the wreath product $(\Z/m\Z) \wr B_n$ \cite[2.1 Surface braid groups]{EsterovLangBraid}. The braid group
$B_n$ acts on $[0;n-1]$ by permutation. This action extends to the group 
$K=\prod_{i=1}^n (\Z/m\Z)_i.$
We define the group $(\Z/m\Z) \wr B_n$ as the semi-direct product $K \rtimes B_n$
with respect to the extended action. Elements in $(\Z/m\Z) \wr B_n$  are braids in $B_{nm}$ each strand of which is decorated with an element in $\Z/m\Z.$ The product of two
elements is performed by first concatenation of two braids, then adding decorations
associated to each strand modulo $m$
The following is a direct consequence of 
Corollary \ref{BnfBn}. 

\begin{corollary}\label{wreathbraid}
The braid group $B_{nm}(f) \leq B_{nm}$ \eqref{exactseq} obtained from roots of the equation \eqref{f(X,Y)1}
is isomorphic to the wreath product
\begin{equation}\label{wreath}
  B_{nm}(f) \simeq (\Z/m\Z) \wr B_n \lneq B_{nm} ,
\end{equation}
while the Galois group $Gal(f)$ of \eqref{f(X,Y)1} admits the corresponding
description thanks to  \eqref{exactseq}
\begin{equation}\label{Galoiswreath}
 Gal(f) \simeq (\Z/m\Z) \wr {\mathfrak S}_n \lneq {\mathfrak S}_{nm}.
\end{equation}
\end{corollary}

\vspace{\fill}

\begin{minipage}[t]{12.2cm}
\begin{flushleft}
{\footnotesize
Susumu TANAB\'E (corresponding author)\\
Department of Discrete Mathematics,\\
Moscow Institute of Physics and Technology,\\
141701, Dolgoprudny, Moscow Region, Russian Federation\\
{\it E-mails}:   {tanabe.s@mipt.ru, tanabesusumu@hotmail.com}\\
ORCID id:  0000-0003-0489-2838  \\

\vspace{0.5cm}
Mutlu KO\c{C}AR\\
Department of Mathematics,\\
Ko\c{c} University,\\
Rumelifeneri Yolu, 34450, Istanbul, Turkey.\\
{\it E-mails}:   {mkocar19@ku.edu.tr, mutlukocar13@hotmail.com}
}\\
\end{flushleft}
\end{minipage}

\end{document}